\documentclass[11pt]{amsart}

\usepackage[latin1]{inputenc}
\usepackage[T1]{fontenc}
\usepackage[english]{babel}

\usepackage{mathrsfs}
\usepackage{amscd}
\usepackage{amsfonts}
\usepackage{amsmath}
\usepackage{amssymb}
\usepackage{amstext}
\usepackage{amsthm}
\usepackage{amsbsy}

\usepackage{cleveref}
\usepackage{tikz-cd}
\usepackage{xspace}
\usepackage[all]{xy}
\usepackage{graphicx}
\usepackage{url}
\usepackage{latexsym}

\newenvironment{claimproof}[1]{\par\noindent\emph{Proof:}\space#1}{\hfill $\square_{Claim}$ \bigskip}

\usepackage{comment}

\newtheorem{theo}{Theorem}[section]
\newtheorem{theorem}[theo]{Theorem}
\newtheorem{lemma}[theo]{Lemma}

\newtheorem{corollary}[theo]{Corollary}
\newtheorem{fact}[theo]{Fact}

\newtheorem{claim}[theo]{Claim}

\theoremstyle{definition}
\newtheorem{dfn}[theo]{Definition}
\newtheorem{remark}[theo]{Remark}

\renewcommand{\leq}{\leqslant}
\renewcommand{\geq}{\geqslant}

\begin{document}

\author{Alf Onshuus}
\address{Departamento de Matem\'aticas, Universidad de los Andes, Bogot\'a, Colombia}
\email{aonshuus@uniandes.edu.co}

\thanks{I would like to give special thanks to Gal Binyamini and Patrick Speissegger for their generous and patient explanations regarding Rolle leaves. I am also very grateful to Annalisa Conversano, Sergei Starchenko, Thomas Scanlon and the anonymous referees for their valuable comments. This research was partially done at MFO during the \emph{Model Theory: Combinatorics, Groups, Valued Fields and Neostability} 2023 workshop. }

\title{Continuous isomorphisms between groups definable in o-minimal expansions of the real field}


\begin{abstract}
In this paper we study the relation between the category of real Lie groups and that of groups definable in o-minimal expansions of the real field, which we will refer to as ``definable groups''. With this terminology, it is known (\cite{Pi88}) that any definable group is a Lie group, and in \cite{COP} a complete characterization of when a Lie group is \emph{Lie isomorphic} to a definable group'' was given. We continue the analysis by explaining when a Lie isomorphism between definable groups is definable.

Among other things, we generalize Wilkie's result on the o-minimality of the exponential function (\cite{Wilkie}) by completely characterizing when, given an o-minimal expansion $\mathcal R$ of the real field  and a Lie isomorphisms $\phi$ between two $\mathcal R$-definable groups $G_1, G_2$,  $\phi$ can be added to the language of $\mathcal R$ preserving o-minimality. We also prove that any definable group $G$ can be endowed with an analytic manifold structure definable in $\mathcal R_{\text{Pfaff}}$ that makes it an analytic group.
 \end{abstract}

\maketitle


 \section{Introduction and preliminaries}
It is known that for any $k\in \mathbb N$, any group definable in an o-minimal expansion of the real field can be definably endowed with a  $C^k$-Lie group structure (see \cite{Pi88} or \cite{PPSlinear}). It is also shown in these papers that definable morphisms between groups definable in such structures are $C^k$-Lie maps. This paper investigates the converse: when is a Lie morphism between definable groups definable.

\medskip

We will fix $k=2$ above, so that by a \emph{Lie group}  we mean a group with a compatible $C^2$-manifold structure and by a \emph{Lie group homomorphism} we mean a $C^2$-map between the manifold structure of two Lie groups. Given a group definable in an o-minimal expansion of the real closed field we will understand it to be a Lie group with the definable $C^2$-Lie structure.

Let $G_1$ and $G_2$ be two groups which are definable in o-minimal expansions of the real field, and let  $\phi$ be a Lie isomorphism between them. If $G_1$ and $G_2$ are abelian, the question of whether or not $\phi$ is a definable isomorphism has been studied for a long time. 

It is known that in the real field with no added structure the exponential map between the additive and multiplicative group is not definable, any one dimensional torsion fee group is definably isomorphic to either the additive or the multiplicative group. Also in the real field, if $A$ is an abelian variety then $A(\mathbb R)$ is a definable but it does not split as products of one dimensional tori (despite being isomorphic as Lie groups). 

In o-minimal structures expanding the real exponential field, much less is known. It is an open question whether or not any two one dimensional torsion free ordered groups are definably isomorphic (see \cite{MiSt} for a survey around this question) and it is also open whether any torsion free abelian two dimensional definable group splits definably as the product of one dimensional subgroups (\cite{PStorsionfree} provides a good vision of the status of the question and gives answers in some cases).

\medskip

Throughout the paper, we fix an o-minimal expansion $\mathcal R$ of the real field. By ``definable group'' and ``definable map'' we mean, respectively, ``$\mathcal R$-definable group'', and ``$\mathcal R$-definable map''.

Given a definable group $G$ we will use the following notations. $Z_G$ will be the center of $G$, $G^0$ the connected component of $G$. The radical $R_G$ of $G$ is the maximal connected solvable subgroup of $G^0$, $U_G$ will denote the maximal simply connected normal subgroup of $R_G$, $K_G:=R_G/U_G$ and $S_G:=G^0/R_G$.

$Z_G, R_G$ and $U_G$ are definable, and $U_G$ is the maximal normal torsion free subgroup of $R_G$ (\cite{edmundo}), $K_G$ is definably compact and abelian and $S_G$ is semisimple with finite center.

\medskip

Our first result is that in any o-minimal expansion of the real field, definability of Lie homomorphisms between definable \emph{abelian} groups is the only obstruction to having a Lie homomorphism be definable:

\begin{theo}\label{Theorem}
Assume that $\mathcal R$ includes the exponential function. Let $\alpha$ be a surjective Lie homomorphism between two definable groups $G_1$ and $G_2$, and assume that
the restriction of $\alpha$ to any definable abelian subgroup of $R_{G_1}$ is definable. Then $\alpha$ is definable.
\end{theo}

We will in fact prove something slightly stronger, since we will only need to assume definability of the restriction of $\alpha$ to central and near central subgroups.

\begin{dfn}
We will say that a definable connected subgroup $A$ of a definable $G$ is \emph{central} if $A\subseteq Z_G$ and that $A$ is \emph{near central} if $A/(Z_G\cap A)$ has dimension less than or equal to one.
\end{dfn}

If $A$ is near central, torsion free, and $A\not\subseteq Z_G$, then $Z_G\cup \{a\}\subseteq C_G(C_G(a))$ for any $a\in A\setminus (A\cap Z_G)$ where $C_G(X):=\{g\in G \mid x\in X\Rightarrow gx=xg \}$. It follows by connectedness, definability and dimension hypothesis that $A\subseteq C_G(C_G(a))$. So any near central torsion free subgroup of $G$ is abelian.

\begin{remark}\label{DCC}
  o-minimal theories have \emph{descending chain condition} (\cite{Pi88}) which states that any infinite chain of definable subgroups is stationary. This implies that $C_G(X)=\bigcap_{x\in X} C_G(x)$ is definable for any set $X$.
\end{remark}

\medskip

For the main applications, we will need to work with the expansion $\mathcal R_{\text{Pfaff}}$ of $\mathcal R$, which includes the Pfaffian closure of a structure. We will refer the reader to \cite{Sp} for the definitions of $\mathcal R_{\text{Pfaff}}$, of Pfaffian equations, and of Pfaffian systems, although we will quote the main result of this paper, which we will need:

\begin{fact}\label{Sp}
Let $\mathcal R$ be an o-minimal expansion of the real field. Then there is an o-minimal
expansion $\mathcal R_{\text{Pfaff}}$ of $\mathcal R$ which is closed under solutions to $C^1$-Pfaffian equations.

In particular, if $U\subseteq \mathbb R$ is a definable open set and $f:U\rightarrow \mathbb R$ is a $C^1$-function satisfying a $C^1$-Pfaffian equation,  then $(\mathcal R, f)$ is o-minimal.
\end{fact}

\medskip

The best known example of Fact \ref{Sp} (which predated it) is Wilkie's result that states that adding exponential map (the continuous isomorphism between $(\mathbb R, +)$ and $(\mathbb R^{>0}, \cdot)$ ) to the real field structure preserves o-minimality. The second part of the paper is devoted to understand when a continuous isomorphism between definable groups can be added to the language preserving o-minimality. The one obstruction that we will have for this comes from the following fact, which is Proposition 2.1 in \cite{CoPi}.

\begin{fact}\label{CoPi}
Let $G$ be a group definable in an o-minimal theory. Then there is a normal definable torsion free group $H$ which contains every normal definable torsion-free subgroup of $G$.
\end{fact}

To understand the obstruction this presents to adding isomorphisms, consider the the group $G:=(\mathbb R, +)\times S_1(\mathbb R)$, with $S_1(\mathbb R)$ the circle group represented as the complex elements of norm one endowed with the complex multiplication, so we think of $G$ as a subgroup of $\mathbb R\times \mathbb C$. Consider the automorphism $\phi$ of $G$ given by $(x,t)\mapsto (x, te^{ix})$. It is not hard to see this is a Lie automorphism of $G$, but $\phi$ can't be definable in \emph{any} o-minimal expansion of the real field: Otherwise in such an expansion the subgroups $\phi\left(\mathbb R\times \{e\}\right)$ and $\mathbb R\times \{e\}$ would be torsion free definable subgroups of $G$ violating Fact \ref{CoPi}.

\medskip

Turns out, this is the only restriction to adding continuous isomorphisms while preserving o-minimality:

\begin{theo}\label{ThmHomomorphism}
Let $\mathcal R$ be an o-minimal expansion of the real field, let $G_1$ and $G_2$ be $\mathcal R$-definable groups with $U_1$ and $U_2$ being the maximal torsion free normal definable subgroups of $G_1$ and $G_2$, respectively. Let $\alpha$ be a Lie isomorphism from $G_1$ to $G_2$ such that $\alpha(U_1)=U_2$. Then $\alpha$ is definable in $\mathcal R_{\text{Pfaff}}$. In particular, $\mathcal R\cup \{\alpha\}$ is o-minimal.
\end{theo}

This, combined with results in \cite{COP} will imply that any $\mathcal R$-definable group can be endowed with an analytic Lie group structure definable in $\mathcal R_{\text{Pfaff}}$.

\medskip

\subsection{Preliminaries}

If $G$ is a connected Lie group, then a \emph{Levi decomposition} of $G$ is a decomposition of the form $G = R_G\cdot L_G$ where $L_G$ is a maximal connected semisimple subgroup, called a \emph{Levi subgroup} of $G$. Equivalently, $L_G$ is a minimal subgroup of $G$ such that the image of $L_G$ under the projection $G\rightarrow G/R_G$ is surjective.

We know the following:

\begin{itemize}
\item $R_G\cap L_G$ is discrete.

\item $L_G$ is a maximal connected semisimple (Lie)-subgroup of $G$, and unique up to conjugation (\cite[Theorem 3.18.13]{Varadarajan}).

\item If $G$ is a definable group, then $L_G$ is a maximal connected ind-definable subgroup of $G$. There are examples of definable groups $G$ where $L_G$ is not definable. (\cite[Theorem 1.1]{CP-Levi}).
\end{itemize}

\medskip

We will need the following results from  \cite{COS} and \cite{COP}:

\begin{fact}\label{supersolvable}
Any torsion free definable group $U$ is supersolvable: This is, there are normal subgroups $U_1, U_2, \dots, U_n$ of $U$ (so $U_i\unlhd U$)  with  $U_n=U$, $U_0=\{e\}$ and such that $U_{i+1}/U_i$ a one dimensional torsion free group.
\end{fact}

\begin{fact}
\label{mainOld}
A real Lie group $G$ is Lie isomorphic to a group definable in an o-minimal expansion of the real field if and only if both $G$ and its center have finitely many connected components and the maximal simply connected normal subgroup $U_G$ of the solvable radical of $G$ is supersolvable.
\end{fact}

\begin{fact}\label{ZS}
If $G$ is definable then $Z_G\cdot L_G$ and $(Z_G)^0 \cdot L_G$ are definable subgroups of $G$.
\end{fact}

We will say that a definable group $G$ is \emph{linear} if it has a definable and faithful finite dimensional representation $\rho$. A \emph{matrix group} will be defined to be a definable subgroup of
$Gl_n(\mathbb R)$ for some $n\in \mathbb N$.

The following is Remark 4.4 in \cite{PPSlinear}.

\begin{fact}\label{PPS}
If $S$ is a semisimple matrix group, then it is semialgebraic.
\end{fact}

We will need some facts about the commutator subgroups.

\begin{dfn}
  For any subset $Y$ of a group $G$, define $[Y,Y]=\bigcup_{n\in \mathbb N}[Y, Y]_n$ where  $[Y, Y]_1=\{xyx^{-1}y^{-1}\mid x,y\in Y\}$ and $[Y, Y]_{n+1}=[Y, Y]_n\cdot [Y, Y]_1$.
  
\noindent  If $G$ is any group, the \emph{commutator subgroup of $G$} is $[G,G]$. 
  
\noindent   A group \emph{perfect} if it is equal to its commutator subgroup.
\end{dfn}

Finally, we will need the following fact, which is Claim $3.1(v)$ in \cite{HPP}.

\begin{fact}\label{HPP3.1}
Let $G$ be a definably connected semisimple group definable in an arbitrary o-minimal structure. Then $G$ is perfect and moreover there is some $r$ such that $[G,G]_r=[G,G]=G$.
\end{fact}

\section{Proof of Theorem \ref{Theorem}}

Let $\alpha: G\rightarrow G'$ be a surjective Lie homomorphism of definable groups such that the restriction of $\alpha$ to any near central (definable) subgroup of $R_G$ is definable. We will prove in this section that $\alpha$ is definable.

\medskip

We will fix $G, G'$ and $\alpha$ as above. We will assume until Subsection \ref{subsection 2.3} that $G=(G, \cdot)$ is connected. Let $L=L_G$ be a Levi subgroup of $G$. To make the notation lighter, let $Z=Z_G$, $Z^0=(Z_G)^0$, $R=R_G$, $U=U_G$, $S=S_G$ and $K=K_G$.

\medskip

To sketch of the proof is as follows: We will prove first that the restriction of $\alpha$ to $U$ is definable, then use this to show that the restriction of $\alpha$ to $R$ is definable. We will then prove that the restriction of $\alpha$ to the definable $Z^0\cdot L$ is definable. Since $G:=L\cdot R=L\cdot Z^0\cdot R$, $\alpha$ is definable whenever its restrictions to $R$ and $L\cdot Z^0$ are, this will complete the proof of the connected case.

\subsection{The restriction of $\alpha$ to $R$ is definable.}

We will need the following.

\begin{dfn}
Let $H$ be a definable subgroup of $G$. We will say that $H$ \emph{definably splits} if there are one dimensional connected definable subgroups $L_1, \dots, L_k$ such that any element of $H$ is a product of elements in $L_i$'s (so that  $H=L_1\cdot L_2\cdot \dots \cdot L_k$).
\end{dfn}

\begin{remark}\label{definable splitting}
Let $H$ be a definable subgroup of $G$ and assume that $H/(Z\cap H)$ definably splits as $H/(Z\cap H)={L_1}\cdot {L_2}\cdot \dots \cdot {L_k}$ (we are abusing notation and using the same notation for the group operations in $G$ and $G/Z$). Let $\pi: G\rightarrow G/Z$ be the natural projection. Then for any $i\leq k$ the group  $\pi^{-1}(L_i)$ is by definition a near central subgroup of $G$ and
\[H=\pi^{-1}\left(H/\left(Z\cap H\right)\right)=\pi^{-1}\left(L_1\right)\cdot \pi^{-1}\left(L_2\right)\cdot \dots \cdot \pi^{-1}\left(L_k\right).\]
\end{remark}

\begin{lemma}\label{restriction to torsion free}
If the restriction of $\alpha$ to any near central subgroup of $U$ is definable, then the restriction of $\alpha$ to $U$ is definable.
\end{lemma}

\begin{proof}
Consider the (faithful) representation of $U/(Z\cap U)$ induced by the adjoint representation of $G$. It is a subgroup of the linear automorphisms of the Lie algebra $\mathfrak g$ of $G$, which is definable by \cite{PPSlinear}. We fix a basis for the vector space $\mathfrak g$ so that $\mathfrak g$ is definably isomorphic to $\mathbb R^n$, and $Ad_G(U)$ is definably isomorphic to its image $Ad^{\mathbb R}_G(U)$ in $Gl_n(\mathbb R)$ induced by this isomorphism. $Ad^{\mathbb R}_G(U)$ is a supersolvable matrix group.

\begin{claim}
Any definable supersolvable matrix group definably splits.
\end{claim}

\begin{claimproof}
Let $H\leq Gl_n(\mathbb R)$ be a supersolvable matrix group. We will prove that $H$ definably splits by induction on $dim(H)$. If $dim(H)=1$ there is nothing to prove.

By supersolvability, let $N$ be $dim(H)-1$ dimensional normal subgroup of $H$.

Let $exp: gl_n(\mathbb R)\rightarrow Gl_n(\mathbb R)$ be the matrix exponentiation, and let $\mathfrak n$ and $\mathfrak h$ be the preimages of $N$ and $H$. Let $\mathfrak l$ be a one dimensional subspace of $\mathfrak h$ such that $\mathfrak h=\mathfrak n+\mathfrak l$.

By Lemma 3.1 in \cite{COS} all eigenvalues of the matrices in $H$ have positive real eigenvalues, which implies (as in Lemma 3.5 of \cite{COS}) that the restriction of $exp$ to $\mathfrak h$ is a definable map (recall that Theorem \ref{Theorem} assumes that $\mathcal R$ includes the exponential function).

By \cite{dixmier} (and supersolvability) $exp$ is a diffeomorphism from $\mathfrak h$ to $H$. So $exp(\mathfrak n)=N$ and $exp(l)$ is a one dimensional torsion free subgroup of $H$ not contained in $N$. This implies that $H=N\cdot exp(l)$. By induction hypothesis $N$ definably splits, so $H$ definably splits, as required.\end{claimproof}

$U/(Z\cap U)$ is definably isomorphic to $Ad^{\mathbb R}_G(U)$ which by the claim definably splits. By Remark \ref{definable splitting}, $U$ is the product of near central subgroups of $G$. By hypothesis, the restriction of $\alpha$ to each of these near central subgroups is definable, which implies that the restriction of $\alpha$ to $U$ is definable.
\end{proof}

\medskip

Recall that $0\rightarrow U\rightarrow R\rightarrow K\rightarrow 0$ is an exact sequence of definable maps, which splits in the Lie category. Unfortunately, $R$ doesn't always have a definable subgroup isomorphic to $K$ (\cite{Strz} provides a very nice example). So we need to work with definable subgroups of $R$ that are more general than definably compact ones. The following are from \cite{Strz}.

\begin{dfn}
Let $E$ be the o-minimal Euler characteristic. A definable subgroup $H$ of $R$ is \emph{a 0-Sylow subgroup} if it is a maximal definable subgroup such that $E(H/H')=0$ for every proper definable subgroup $H'$ of $H$. \end{dfn}

\begin{fact}
0-Sylow subgroups of definable groups are abelian.
\end{fact}

The following is Proposition 3.1 in \cite{ConvMax}):

\begin{fact}\label{Conversano T}
Let $T$ be a 0-Sylow subgroup of the definably connected solvable group $R$, and let $U$ be the maximal torsion free subgroup of $R$. Then $R=T\cdot U$.
\end{fact}

\begin{lemma}\label{restriction to R}
If the restriction of $\alpha$ to any abelian subgroup of $G$ is definable, the restriction of $\alpha$ to $R$ is definable.
\end{lemma}

\begin{proof}
Let $T$ be a 0-Sylow subgroup of $R$. As before, the group $T/(Z\cap T)$ is definable, definably connected, and abelian. As in the proof of Lemma \ref{restriction to torsion free} it
is definably isomorphic to a matrix group, so by  Fact 3.1 and Proposition 3.8 in \cite{PPSlinear} it definably splits. By Remark \ref{definable splitting} the group $T$ is a product of near central subgroups of $G$ which by hypothesis on $\alpha$ implies that the restriction of $\alpha$ to $T$ is definable. So $R=T\cdot U$ and the restrictions of $\alpha$ to both $U$ (Lemma \ref{restriction to torsion free}) and $T$ are definable, the restriction of $\alpha$ to $R$ is definable.
\end{proof}

\subsection{Definability of the restriction of $\alpha$ to $L\cdot Z^0$.}

\begin{lemma}\label{H1}
If the restriction of $\alpha$ to $Z^0$ is definable, then the restriction of $\alpha$ to $L\cdot Z^0$ is definable.
\end{lemma}

\begin{proof}
Let $H:=L\cdot Z^0$ which, by Fact \ref{ZS}, is definable. Let $H'=\alpha(H)=\alpha(L)\cdot \alpha(Z^0)$. Let $\pi$ and $\pi'$  be the projection maps from $H$ and $H'$  to $S_H:=G/R$ and $S_H':=\alpha(G')/\alpha(R_{G'})$, respectively. Let $\sigma$ be the homomorphism between $S_H$ and $S'_H$ induced by $\alpha$. By construction, $S_H$ is semisimple.  

\begin{claim}
$S_H'=\alpha(H)/\alpha(Z(H))$.
\end{claim}

\begin{claimproof}
  $\pi(H)=S_H$ and $\sigma$ is surjective, so $\pi'(\alpha(H'))=S_H'$.

  $S_H'$ is the image of the semisimple group $S_H$, so  $S_H'$ is semisimple and we have $\alpha(Z(H))=\alpha(\ker(\pi_1))\geq R_{H'}$. On the other hand, $\alpha(Z(H))\subseteq Z(H')$ and since $Z(H)=Z^0$, $Z(H')$ is connected so $Z(H')\leq R_{H'}$. So $\alpha(Z(H))=Z(H')$.

  By construction $S_H= H/Z(H)$ so $S'(H)=\sigma(S_H)=\alpha(H)/\alpha(Z(H))$, the result follows.
\end{claimproof}

By Lemma \ref{restriction to R} (and definable choice) $S_H'$ is definable, so $H'$ is a central extension of the definable $S_H'$ by $\alpha(Z^0)$, which is definable by hypothesis. Fact \ref{ZS} implies that $H'$ is definable.

\medskip

As before, the images of $S_H$ and $S'_H$ in the adjoint representations of $H$ and $H'$ are faithful, which implies that both groups are definably isomorphic to matrix groups. By Fact \ref{PPS} the graph of the induced map between these matrix copies is semialgebraic, so that  $\sigma$ is definable. By definition $\sigma(S_H)=S'_H$ so $\sigma$ is a definable continuous quotient map.

\medskip

Fact \ref{HPP3.1} implies that $S_H=[S_H,S_H]=[S_H, S_H]_k$  for fixed some $k\in \mathbb N$. $L$ is a Levi subgroup of $H$ which is a central extension of a semisimple group so $[H,H]=L$. By definition $\pi$ maps $[H, H]_i$ surjectively onto $[S_H, S_H]_i$ for all $i\in \mathbb N$, so $[H,H]_k$ is a definable subset of $L$ which is projected by $\pi$ onto $S_L$. 

By definable choice, we can take a definable global section $\rho:S_H\rightarrow [H,H]_k$ of $\pi$ so that in particular the image of $\rho$ is contained in $L$. Similarly, we can find a definable global section $\rho':S_H'\rightarrow H'$ of $\pi'$ such that $\rho'$ is contained in $\alpha(L)=[H',H']$.

We have the following diagram:
\[
\begin{tikzcd}
H \arrow[r, "\alpha"] \arrow[d, "\pi"]
& H' \arrow[d, "\pi'" ] \\
S_H \arrow[u, bend left, "\rho"] \arrow[r, "\sigma"]
&  S_H' \arrow[u, bend left, "\rho'"].
\end{tikzcd}
\]

By cell decomposition, there are open dense $V\subseteq S_H$ and $V'\subseteq S_H'$ in which $\rho$ and $\rho'$ are continuous. $\sigma$ is a quotient map (open and continuous) so $\sigma^{-1}(V')$ and $\sigma(V)$ are open dense subsets of $S_H$ and $S_H'$, respectively. Replacing $V$ with $V\cap \sigma^{-1}(V')$ and $V'$ with $V'\cap \sigma(V')$ we can assume that $\sigma$ maps $V$ onto $V'$.

Consider the map $\phi$ from $V$ to $H'$ given by $c\mapsto (\rho'( \sigma(c)))\cdot (\alpha( \rho(c)))^{-1}.$

Since $\pi'(\rho'( \sigma(c)))=\sigma(c)= \pi'(\alpha(\rho(c)))$, we have
\[(\rho'(  \sigma(c)))\cdot (\alpha(\rho(c)))^{-1}\in Ker(\pi')=Z(\alpha(H)).\] The images of both $\rho'(\sigma(c))$ and $\alpha(\rho(c))$ are contained in $\alpha(L)$, which implies that $\phi$ is continuous with image contained in the discrete set $\alpha(L)\cap Z(H')$.

The set
\[\{\phi^{-1}(x)\}_{x\in \alpha(L)\cap Z(H')}\] is therefore an open cover of $V$ by disjoint open sets. It follows that $\phi$ is constant in each of the finitely many (definable) connected components of $V$, so  $\phi$ is a definable map with finite image.

Restricted to $V$ we have then $\alpha\circ \rho= (\rho'\circ \sigma)\cdot \phi$ so that $\alpha\circ \rho$ is a definable map from $V$ to $H'$ and the restriction of $\alpha$ to $\rho(V)$ is definable. Using multiplication in $H$ we know that the restriction of $\alpha$ to $L_0:=\rho(V)\cdot \rho(V)$ is definable. The projection of $L_0$ to $S_H$ contains $V\cdot V$ which, since $V$ was open dense, is all of $S_H$. It follows that $L_0\cdot Z^0=L\cdot Z^0$ and $\alpha$ is definable, as required.
\end{proof}

\medskip

\subsection{Proof of Theorem \ref{Theorem}}\label{subsection 2.3}

\begin{proof}[Proof of Theorem \ref{Theorem}]
If $G$ is connected, then $G=L\cdot R=(L\cdot Z^0)\cdot R$. The restrictions of $\alpha$ to $L\cdot Z^0$ and $R$ are both definable, so $\alpha$ is definable.

\medskip

For the general case, $G/G^0$ is finite so there is a finite subset $F$ of $G$ such that $G=F\cdot G^0$. Since the restriction of $\alpha$ to $G^0$ is definable, and the restriction of $\alpha$ to $F$ is clearly definable, we get that $\alpha$ is definable as required.
\end{proof}

\section{Adding the Lie isomorphisms between definable groups in expansions of the real field.}\label{Adding}

In this section we characterize when, given any o-minimal expansion $\mathcal R$ of the field of real numbers, a Lie isomorphism $\phi$ between definable groups $G_1$ and $G_2$ can be added preserving o-minimality.

As mentioned in the introduction, if $\mathcal R\cup \{\phi\}$ is o-minimal, and $U_1$ and $U_2$ are the maximal torsion free normal subgroups of $G_1$ and $G_2$, then $\phi(U_1)=U_2$ since otherwise the subgroups $U_2$ and $\phi(U_1)$ would be definable in $\mathcal R\cup \{\phi\}$ contradicting Fact \ref{CoPi}.  

\medskip

We will prove a converse of this observation. If $\phi$ is such that the image of $U_1$ is $U_2$, then $\phi$ is definable in $\mathcal R_{\text{Pfaff}}$, and by Fact \ref{Sp} $(\mathcal R, \phi)$ is o-minimal.

\medskip

The following is a very useful result, which appears as a lemma in (\cite{PSS}).

\begin{fact}\label{PSS1}
Let $\mathcal R$ be an o-minimal expansion of the field of real numbers, and let
$(A, <, U, *)$ be a group interval of class $C^1$ definable in $\mathcal R$. Assume that $(B, <, V, +)$
is a sub-group interval of $(\mathbb R, <,\mathbb R^2, +)$ and that $f : (B, <, V, +)\rightarrow (A, <, U, *)$
is a $C^1$-isomorphism. Then the structure $(\mathcal R, f)$ is o-minimal.
\end{fact}

The proof in \cite{PSS} shows that $f$ satisfies the Pfaffian system on $V$, given by $f'=h(f(x))$ where $h(a) := \frac{\partial (a*b)}{\partial b}\mid_{(a, 1)}\cdot f'(0)$. By Fact \ref{Sp}, $(\mathcal R, f)$ is o-minimal.

We abuse notation and by $(\mathbb R/\mathbb Z, +)$ we denote the definable group with universe $(-1/2, 1/2]$ and group operation is
\[x\oplus y= \begin{cases}
      x+y-1 & x+y> 1/2 \\
      x+y & -1/2< x+y\leq  1/2 \\
      x+y+1 &  x\leq -1/2
   \end{cases}
\]
The following follows from the above fact.

\begin{fact}\label{PSS}
Let $\mathcal R$ be an o-minimal expansion of the real field which is closed under solutions of Pfaffian equations. Then the following hold:

\begin{itemize}
\item Let $(G, *)$ be a torsion free one dimensional definable group. Then $(G, *)$ is $\mathcal R$-definably isomorphic to $(\mathbb R, +)$.

\item Let $(G, *)$ be a compact one dimensional definable group. Then $(G, *)$ is $\mathcal R$-definably isomorphic to $(\mathbb R/\mathbb Z, +)$.
\end{itemize}
\end{fact}

\begin{proof}
Since ordered groups are particular cases of group intervals, the first item follows straight from Fact \ref{PSS1}. For the second one, recall that by  \cite{PPSlinear} in any group $G$ definable  in an o-minimal expansion of the real field we can define a $C^1$-structure that makes $G$ a Lie group and that over the real field definable compactness implies compactness. In particular, if $G$ is a definably compact abelian one dimensional Lie group and the Lie exponential map $Exp$ from the Lie algebra $(\mathbb R,+)$ to $G$ is a $C^1$ surjective group morphism. The kernel must be  discrete and generated by its first positive element. After a linear transformation, we may assume that the kernel is $\mathbb Z$. The map $Exp$ restricted to the interval $(-1/2, 1/2)$ (with the standard order inherited from $\mathbb R$) will be a $C^1$ group interval isomorphism, which by Fact \ref{PSS1} is definable. This map can be extended definably to a group isomorphism from $\mathbb R/\mathbb Z$ to $G$.
\end{proof}

\subsection{Lie isomorphisms between definable torsion free abelian subgroups} \label{ss abelian torsion free}

\begin{lemma}\label{splitting abelian 2}
Let $(H, \odot_H)$ be a two dimensional abelian torsion free group definable in an o-minimal expansion $\mathcal R$ of  the real closed field, and let $U$ be a definable one dimensional subgroup of $H$. Then the sequence
\[ 0\rightarrow U \rightarrow H\rightarrow H/U\rightarrow 0\] definably splits in  $\mathcal R_{\text{Pfaff}}$. This is, there is a definable subgroup $H'$ of $H$ which projects bijectively onto $H/U$ (and such that $H=U\odot H'$).
\end{lemma}

\begin{proof}
By  Fact \ref{PSS} in $\mathcal R_{\text{Pfaff}}$ both $H$ and $H/U$ are $\mathcal R_{\text{Pfaff}}$-definably  isomorphic to $(\mathbb R, +)$. By definable choice,  we may assume that the universe of $H$ is $\mathbb  R \times \mathbb R$ and that the identity of $H$ is $(0,0)$. The group structure will be given by
\[(a,b)\odot (c,d)=\left(\left(a+c\right), \left(b+d+F\left(a,b\right)\right)\right)\] for any $a,b,c,d\in \mathbb R$, where $F(x,y)$ is a 2-cocycle. Once the universe is fixed, $F$ can be computed directly from the group operation, so it is definable in $\mathcal R_{\text{Pfaff}}$.

In the category of Lie groups all abelian groups split, so there is a Lie group section, 
i.e. a  $C^1$-function $\gamma: \mathbb R\rightarrow \mathbb R$ such that $(x,\gamma(x))$ is a subgroup of $H$. So $(x,\gamma(x))\odot (y, \gamma(y))=(x+y, \gamma(x+y))$, so \[F(x,y)=\gamma (x+y)-\gamma(x)-\gamma(y).\]

Then
\[
\gamma'(x)= \lim_{h \to 0} \frac {\gamma(x+h)-\gamma(x)}{h}= \lim_{h \to 0} \frac {F(x,h)}{h}+ \frac {\gamma(h)}{h}.\]

Since $\lim_{x \to 0} \frac {\gamma(h)}{h}=\gamma'(0)$ is a constant, $\lim_{x \to 0} \frac {F(x,h)}{h}$ exists and $\gamma'(x)$ is definable. So $\gamma(x)$ is definable in  $\mathcal R_{\text{Pfaff}}$, which implies that $H$ definably splits witnessed by  $H':=\{(x,\gamma(x)\mid x\in \mathbb R)\}$.
\end{proof}

\begin{corollary}\label{splitting abelian}
Let $(H, \odot_H)$ be a torsion free abelian group definable in an o-minimal expansion  $\mathcal R$ of  the real closed field. Then $H$ definably splits in $\mathcal R_{\text{Pfaff}}$.

In particular, $H$ is $\mathcal R_{\text{Pfaff}}$-definably isomorphic to $(\mathbb R^n, +)$ for some $n$ and every Lie subgroup of $H$ is definable.
\end{corollary}

\begin{proof}
By \cite{PSOneDim}, any torsion free definable group (in fact any non compact definable group) contains a one dimensional definable subgroup, so Lemma \ref{splitting abelian 2} implies that any definable 2-dimensional abelian torsion free group splits in $\mathcal R_{\text{Pfaff}}$. We also know that in $\mathcal R_{\text{Pfaff}}$ any abelian torsion free group is definably isomorphic to $(\mathbb R, +)$.

The rest is a simple induction argument (which can be found in \cite{PStorsionfree}). Assume we have the result for $n\geq i\geq 2$ and let $H$ have dimension $n+1$. By \cite{PSOneDim} there is a definable one dimensional subgroup $U$ of $H$. Now, $H/U$ splits into products of one dimensional groups $R_1\times\dots \times R_n$. The preimage $H_i$ of $R_i$ in $H$ has dimension 2, so by Lemma \ref{splitting abelian 2} it splits so that there is a definable isomorphism $\phi_i$ from $R_i\times U$ to $H_i$. The map $\phi: R_1\times \dots \times R_n\times U$ defined by
\[
\phi(x_1, \dots, x_n, g)\mapsto \phi_1(x_1,g)\odot_H\dots \odot_H \phi_1(x_n,g)
\]
is therefore a group isomorphism.

\medskip

Since $H$ is a direct product of definable torsion free one dimensional subgroups and by Fact \ref{PSS} each one is $\mathcal R_{\text{Pfaff}}$-definably isomorphic to $(\mathbb R,+)$, it follows that $H$ is $\mathcal R_{\text{Pfaff}}$-definably isomorphic to $(\mathbb R^{n+1}, +)$. It is therefore a vector space and any Lie subgroup is a vector subspace, so it is definable. This completes the proof of the corollary.
\end{proof}

\subsection{Compact abelian subgroups} \label{ss compact}

We will fix $(A, \odot_A)$ be a to be a $\mathcal R$-definable abelian group. Let $\mathfrak a$ be its Lie algebra (see \cite{PPSlinear} for an introduction to Lie algebras in the o-minimal context). The results in this subsection will be corollaries to the following proposition, which was pointed to us by Gal Binyamini.

\begin{lemma}\label{Binyamini}
Let $H$ be a connected Lie subgroup of $A$. Then there is a neighborhood $W$ of the identity $e_A$ of $A$ such that $H\cap W$ is definable in $\mathcal R_{\text{Pfaff}}$.
\end{lemma}

\begin{proof}
Let $n$ be the dimension of $A$ and let $\mathfrak a$ and $\mathfrak h$ be the Lie algebras of $A$ and $H$, respectively. Both $A/H$ and $H$ are abelian Lie groups, so by \cite[Theorem 2.12]{OlVi93} they split as a product of a cartesian power of $\mathbb R$ with a cartesian power of the circle group $S_1(\mathbb R)$. We can therefore find a sequence of Lie subgroups
\[A_n=\{e\}\leq A_{n-1}\leq \dots \leq A_l\leq A_{l+1}\leq\dots \leq A_0=A\] of $A$ where $A_l=H$ and $A_i/A_{i+1}$ is a one dimensional group isomorphic to either $(\mathbb R,+)$ or to $S_1(\mathbb R)$.

Let $\mathfrak a_l$ be the Lie algebra of $A_l$, and let $v_1, \dots, v_n\in \mathfrak a$ be elements in $\mathfrak{a}_0$ be such that $\mathfrak{a}_j=Span\left(\left\{v_1, \dots, v_{n-j}\right\}\right)$, so that in particular $\mathfrak h:=\mathfrak{a}_k=Span\left(\left\{v_1, \dots, v_{n-k}\right\}\right)$ is the Lie algebra of $K$.

In the context of \cite[Chapter 10]{Spivak}, let $T(A_l)$ be the tangent bundle of $A_l$, and let $d_l:A\rightarrow T(A)$ be  distribution  associated with the Lie algebra $\mathfrak{a}_l$, so that $d_l(x)=L_x(\mathfrak{a}_l)$ where $L_x$ is the derivative of the left translation (in $A_l$) by the element $x$.

Now, $d:=(d_0, d_1, \dots, d_n)$ is a nested tuple of distributions of $A$, each of which is integrable witnessed by the nested integral manifold $\mathcal A:=(A_0, A_1, \dots, A_n)$ of $d$. So $\mathcal A$ is a nested leaf of $d$.

We now move to the context of \cite{LiSp}.

\begin{claim}
For some neighborhood $W$ of $e_A$, $\mathcal A\cap W=(A_0\cap W, A_1\cap W, \dots, A_n\cap W)$ is a nested Rolle leaf of the restriction of $d$ to $W$.
\end{claim}

\begin{claimproof}
Recall that the nested leaf $(A_0\cap W, A_1\cap W, \dots, A_n\cap W)$ of the nested distribution $d|_W$ is Rolle if it satisfies the separating condition: For each $l$ we have that $A_{l+1}\cap W$ is a closed submanifold of $A_l\cap W$ and for any $C^1$-curve $\gamma:[0,1]\mapsto A_l\cap W$ with $\gamma(0), \gamma(1) \in A_{l+1}$ there exists a $t\in [0,1]$ such that
$\gamma'(t)\in d_{l}(\gamma(t))$.

Now, by hypothesis $\mathfrak a$ is abelian, so it is isomorphic to $\mathbb R^n$, which implies (\cite{Spivak} Theorem 5) there is a local $C^2$-diffeomorphism between $A$ and $\mathfrak a$. This is, there is a neighborhood $W$ of $A$ and a diffeomorphism sending $A_l\cap W$ into $\mathfrak a_l$. If we send the basis $v_1, \dots v_n$ to the canonical basis $e_1\dots e_n$ of $\mathbb R^n$, the diffeomorphisms sends $d$ to the nested trivial distribution  over $\mathbb R^n$ whose nested integral manifold $\langle \mathbb R^{n-i}\rangle$ (and the restriction to the image of $W$) is trivially a nested Rolle leaf. The separating condition is, by definition, preserved under diffeomorphisms, so $(A_0\cap W, A_1\cap W, \dots, A_n\cap W)$ is a nested Rolle leaf of the restriction of $d$ to $W$.
\end{claimproof}

By \cite{PPSlinear} the Lie algebra $\mathfrak{a}$ is $\mathcal R$-definable, and subspaces of vector spaces are always definable, which implies that $\mathfrak{a}_l$ is $\mathcal R$-definable for all $l$. Right translation by an element $x\in A$ is of course $\mathcal R$-definable, and therefore so is its differential $L_x$, which implies that the distribution $d$ is $\mathcal R$-definable. By Proposition 3.6 in \cite{LiSp} the definability of the Rolle leaf $A_l\cap W$ in an o-minimal expansion $\mathcal R'$ of the real field implies that $A_{l+1}\cap W$ is definable in $\mathcal R'_{\text{Pfaff}}$. By induction $A_l\cap W$ is definable in $\mathcal R_{\text{Pfaff}}$ for all $l$, so in particular $H\cap W=A_k\cap W$ is $\mathcal R_{\text{Pfaff}}$-definable, as required.
\end{proof}

\begin{corollary}\label{compact}
The following hold.

\begin{enumerate}
\item Any compact connected subgroup of $A$ is $\mathcal R_{\text{Pfaff}}$-definable.

\item If $B$ is any $\mathcal R$-definable compact abelian group and $\phi:B\rightarrow A$ is a Lie isomorphism, then the restriction of $\phi$ to some open neighborhood of the identity is definable in $\mathcal R_{\text{Pfaff}}$.
\end{enumerate}
\end{corollary}

\begin{proof}
Let $K$ be any compact subgroup of $A$. By compactness the connected component of $K$ has finite index and the definability of $K$ is equivalent to the definability of its connected component, so we may assume that $K$ is connected. By Lemma \ref{Binyamini} for some open $W$ the set $K\cap W$ is definable in $\mathcal R_{\text{Pfaff}}$. By compactness and connectedness,
\[
K=\underbrace{(K\cap W)\odot_A \dots \odot_A (K\cap W)}_{m-\text{times}}
\] for some $m\in \mathbb N$. Since multiplication in $A$ is definable, (1) follows.

For the second item, notice that the graph of $\phi$ is a compact Lie subgroup of the $\mathcal R$-definable abelian $B\times A$.
\end{proof}

\subsection{Proof of Theorem \ref{ThmHomomorphism}}

We now prove Theorem \ref{ThmHomomorphism}:

\begin{theorem}\label{ThmAdding}
Let $G_1$ and $G_2$ be groups definable in $\mathcal R$ with $U_1$ and $U_2$ their maximal normal torsion free subgroups, and let $\alpha: G_1\rightarrow G_2$ be a Lie isomorphism with $\alpha(U_1)=U_2$. Then
 $\alpha$ is definable in $\mathcal R_{\text{Pfaff}}$. In particular $\mathcal R\cup\{\alpha\}$ is o-minimal.
\end{theorem}

\begin{proof}
Let $R$ be the solvable radical of $G_1$, $Z^0$ the connected component of the center of $G_1$ and $L$ be a Levi subgroup of $G_1$. Since $G_1=(R)\cdot (L\cdot Z^0)$ and both $R$ and   $L\cdot Z^0$ are definable, it is enough to prove that the restriction of $\alpha$ to $R$ and to $L\cdot Z^0$ are definable in  $\mathcal R_{\text{Pfaff}}$.

We begin by proving that the restriction of $\alpha$ to $R$ is definable in  $\mathcal R_{\text{Pfaff}}$. 

\begin{claim}
The restriction of $\alpha$ to $R$ is definable in  $\mathcal R_{\text{Pfaff}}$.
\end{claim}

\begin{claimproof}
Let $K$ be a maximal compact Lie subgroup of $R$. By Corollary \ref{compact} $K$, $\alpha(K)$ and the graph of the restriction of $\alpha$ to $K$ are definable in  $\mathcal R_{\text{Pfaff}}$. Since $R=K\cdot U_1$, we need to prove that the restriction of $\alpha$ to $U_1$ is definable in $\mathcal R_{\text{Pfaff}}$. By Lemma \ref{restriction to torsion free}, it is enough to prove that the restriction of $\alpha$ to any near central subgroup of $U_1$ is definable in $\mathcal R_{\text{Pfaff}}$. Let $A$ be a near central subgroup of $U_1$. Let $a\in A$ be such that $C_{G_1}(C_{G_1}(a))\supseteq A$. Since $\alpha$ is an isomorphism, then $C_{G_2}(C_{G_2}(\alpha(a)))\supseteq \alpha(A)$. This, together with the assumption that $\alpha(U_1)=U_2$, imply that both $A$ and $\alpha(A)$ are subgroups of torsion free abelian definable subgroups of $G_1$ and $G_2$ (namely $C_{G_1}(C_{G_1}(a))\cap U_1$ and $C_{G_2}(C_{G_2}(\alpha(a)))\cap \alpha(U_1)$) and by Corollary \ref{splitting abelian} both $A$ and $\alpha(A)$ are definable and definably isomorphic (in  $\mathcal R_{\text{Pfaff}}$) to $(\mathbb R^k, +)$ for some $k$. This implies that  $\alpha$ is interdefinable in  $\mathcal R_{\text{Pfaff}}$ with a linear automorphism of $(\mathbb R^k, +)$, which is of course definable.
\end{claimproof}

Since $Z^0$ is a subgroup of $R$, the above proof shows that the restriction of $\alpha$ to $Z^0$ is definable in $\mathcal R_{\text{Pfaff}}$,  and Lemma \ref{H1} the implies that  $L\cdot Z^0$ is definable in  $\mathcal R_{\text{Pfaff}}$. This completes the proof of the theorem.\end{proof}

\begin{corollary}\label{analytic}
Let $G$ be a group definable in $\mathcal R$. Then $G$ is $\mathcal R_{\text{Pfaff}}$-definably isomorphic to an analytic group.
\end{corollary}

\begin{proof}
The main result in \cite{COP} implies that $G$ is Lie-isomorphic to a group $G_{\text{def}}$ definable in $\mathbb R_{\text{exp}}$. In $\mathbb R_{\text{exp}}$ all functions are locally analytic. By results in \cite{Pi88} and \cite{PPSlinear} $G_{\text{def}}$ can be definably endowed with a manifold structure that make $G_{\text{def}}$ into an analytic Lie group.

By Theorem \ref{ThmHomomorphism} $G$ is  $\mathcal R_{\text{Pfaff}}$-definably isomorphic to $G_{\text{def}}$.
\end{proof}

\section{Further considerations. sharp-o-minimality}

Lately there have been a lot of research about finding effective bounds for finiteness results in o-minimal structures (see \cite{BJST} and \cite{BNZ}). For this results to be valid in our context (when adding Lie-homomorphisms between definable groups) it might be useful to understand precisely the Pfaffian solutions we are adding. These come in three kinds:

\begin{itemize}
\item We add predicates of Lie compact subgroups of abelian definable groups and function symbols for Lie isomorphisms between definable compact groups. These are definable with restricted solutions of Pfaffian systems, which are known to be sharply-o-minimal.

\item We add isomorphisms between any two definable torsion free one dimensional groups. These include the exponential (if it was not originally definable in $\mathcal R$), but we may need to add more functions (the existence of such groups is still conjectural). The isomorphism between any such group (if it exist) and the additive group will satisfy a Pfaffian equation quite analogous to the one satisfied by the exponential function (see \cite{PSS}). We call these \emph{generalized exponentials.}

\item In order to have splitting for all torsion free definable abelian groups, we needed (in our proof) to add antiderivatives of definable functions.
\end{itemize}

So all our results will hold if we take an o-minimal structure and take the closure with respect to restricted solutions of Pfaffian systems, generalized exponentials, and antiderivatives. We believe that this setting fits in the framework of \cite{BJST}, so that, for example, effective Pila-Wilkie and possibly Wilkie's Conjecture should hold after adding all Lie isomorphisms between definable groups.

\bibliographystyle{abbrv}
\bibliography{DefinableOminimalGroups}
\end{document}